\documentclass[proc]{edpsmath}
\usepackage{amsmath}
\usepackage{amsthm}

\theoremstyle{plain}
\newtheorem{theorem}{Theorem}

\DeclareMathOperator*{\esup}{ess\,sup}

\newcommand{\N}{\mathbb N}
\newcommand{\R}{\mathbb R}
\newcommand{\real}{\mathbb R}
\newcommand{\normal}{\mathcal N}
\newcommand{\ind}{\perp\!\!\!\!\perp}
\newcommand{\Yd}{Y^{\dagger}}
\newcommand{\yd}{y^{\dagger}}
\newcommand{\mmu}{\mathfrak{\mu}}
\newcommand{\mnu}{\mathfrak{\nu}}
\newcommand{\op}[1]{\mathsf #1}
\renewcommand{\P}{\op{P}}
\newcommand{\B}{\op{B}}
\newcommand{\G}{\op{G}}
\newcommand{\Q}{\op{Q}}
\newcommand{\T}{\op{T}}
\newcommand{\ppi}{\mathfrak{\pi}}
\newcommand{\hmu}{\widehat \mu}
\renewcommand{\d}{\mathrm d}
\newcommand{\placeholder}{\, \cdot \,}
\newcommand{\pij}{\pi^{{\rm EK}, J}}
\newcommand{\muj}{\mu^{{\rm EK}, J}}
\newcommand{\muk}{\mu^{\rm EK}}
\newcommand{\nuk}{\pi^{\rm EK}}
\newcommand{\mapT}{\mathfrak{T}}
\newcommand{\tvdg}{\mathrm d_g}
\newcommand{\pp}{\mathfrak{p}}
\newcommand{\range}[2]{\{#1, \dotsc, #2\}}
\newcommand{\winf}[1]{\lVert #1 \rVert_{L^{\infty}_{w}}}
\newcommand{\lipnorm}[1]{| #1 |_{C^{0,1}}}

\begin{document}

\vspace{-3cm}
\title{Statistical accuracy of the ensemble Kalman filter in the near-linear setting}

\author{E.~Calvello}
\address{Department of Computing and Mathematical Sciences, Caltech, USA}

\author{J.~A. Carrillo}
\address{Mathematical Institute, University of Oxford, UK}

\author{F.~Hoffmann}
\address{Department of Computing and Mathematical Sciences, Caltech, USA}

\author{P.~Monmarché}
\address{Laboratoire Jacques-Louis Lions, Sorbonne Universit\'e \& Institut Universitaire de France, France}

\author{A.~M. Stuart}
\address{Department of Computing and Mathematical Sciences, Caltech, USA}

\author{U.~Vaes}
\address{MATHERIALS project-team, Inria Paris \& École des Ponts, IP Paris}

% \thanks{%
%     UV is grateful to Edoardo Calvello and Andrew Stuart for very useful feedback.
%     UV is partially supported by the European Research Council (ERC)
%     under the EU Horizon 2020 programme (grant agreement No 810367),
%     and by the Agence Nationale de la Recherche under grants
%     ANR-21-CE40-0006 (SINEQ) and ANR-23-CE40-0027 (IPSO).
% }

\begin{abstract}
    Estimating the state of a dynamical system from partial and noisy observations is a ubiquitous problem in a large number of applications,
    such as probabilistic weather forecasting and prediction of epidemics.
    Particle filters are a widely adopted approach to the problem
    and provide provably accurate approximations of the statistics of the state,
    but they perform poorly in high dimensions because of weight collapse.
    The ensemble Kalman filter does not suffer from this issue, as it relies on an interacting particle system with equal weights.
    Despite its wide adoption in the geophysical sciences,
    mathematical analysis of the accuracy of this filter is predominantly confined to the setting of linear dynamical models and linear observations operators,
    and analysis beyond the linear Gaussian setting is still in its infancy.
    In this short note,
    we provide an accessible overview of recent work in which the authors take first steps
    to analyze the accuracy of the filter beyond the linear Gaussian setting~\cite{carrillo2022ensemble,carrillo2024statisticalaccuracyapproximatefiltering,calvello2024accuracyensemblekalmanfilter}.
\end{abstract}

\begin{resume}
    Estimer l'état d'un système dynamique à partir d'observations partielles et bruitées est un problème omniprésent dans de nombreuses applications, telles que la prévision météo et la prédiction des épidémies.
    Les filtres particulaires sont largement utilisés pour résoudre ce problème et fournissent des approximations des statistiques de l'état dont la précision peut être démontrée rigoureusement.
    Cependant, ils fonctionnent mal en grande dimension en raison de la disparité des poids associés aux particules.
    Le filtre de Kalman d'ensemble ne souffre pas de ce problème, car il repose sur un système de particules en interaction avec des poids égaux.
    Malgré sa popularité dans les sciences géophysiques, l'analyse mathématique de la précision de ce filtre est principalement limitée au cadre des modèles dynamiques linéaires et des opérateurs d'observation linéaires, et l'analyse au-delà de ce cadre en est encore à ses débuts.
    Dans cette courte note, nous présentons un aperçu accessible de travaux récents dans lesquels
    la précision du filtre est étudiée au-delà du cadre linéaire gaussien~\cite{carrillo2022ensemble,carrillo2024statisticalaccuracyapproximatefiltering,calvello2024accuracyensemblekalmanfilter}.
\end{resume}

\maketitle
\section*{Introduction}
This note is partly based on~\cite{carrillo2024statisticalaccuracyapproximatefiltering}.
Consider a discrete-time dynamical system with state $\{v_n\}_{n \in \N}$ evolving in~$\real^d$,
observed through partial and noisy data~$\{y_n\}_{n \in \mathbb{N}}$ in $\real^K$.
Suppose that the state and data are determined by the following system:
\begin{alignat*}{3}
    \text{ \bf State evolution: } \qquad
    v_{n+1} &= \Psi(v_{n}) + \xi_{n}\, ,\\
    \text{ \bf Data aquisition: } \qquad
    y_{n+1} &= h(v_{n+1}) + \eta_{n+1}\,.
\end{alignat*}
This holds for $n \in \N = \{0,1,2,\cdots\}$.
We assume that the initial state of the system is a Gaussian random variable~$v_0 \sim \normal(m_0, C_0)$
with known mean and covariance.
In addition,
we assume that the noise entering in the state evolution is distributed according to the normal distribution $\xi_{n} \sim \normal(0, \Sigma)$,
and that the noise affecting the observation satisfies $\eta_{n+1} \sim \normal(0, \Gamma)$.
We assume that both noise sequences are i.i.d.\ and
satisfy the following independence assumption:
\[
    v_0 \ind \{\xi_n\}_{n \in \N} \ind \{\eta_{n+1}\}_{n \in \N}.
\]
The objective of the filtering problem is to determine,
and update sequentially as new data arrives,
the probability distribution on state $v_n$ given all the data observed up to time~$n$.
To state the problem more precisely,
we consider a given realization of the data denoted by a dagger and define
\[
    \Yd_n  =\{\yd_\ell\}_{\ell=1}^n, \qquad v_n|\Yd_n \sim \mmu_n.
\]
The objective of filtering methods can then be reformulated precisely as follows:
it is to approximate the \emph{filtering distribution} or \emph{true filter}~$\mmu_n$ and update it sequentially in $n$.

\section{Evolution of the filtering distribution}%
\label{sec:Evolution of the filtering distribution}
In this section,
we show that the law $\mmu_n$ evolves according to a nonlinear and nonautonomous dynamical system on the space of probability measures.
To determine this dynamical system, and the approximation of it by the ensemble Kalman filter discussed in Section~\ref{section:enkf},
we denote by $\mathcal P(\real^r)$ the set all probability measures on~$\real^r$,
and by $\mathcal G(\real^r)$ all the set Gaussian probability measures on $\real^r$.
For simplicity we will use the same symbol for measures and their
densities throughout this article.
To write the evolution of~$\mu_n$ in a compact manner,
it is useful to introduce key operators on probability measures:
linear operator $\P\colon \mathcal P(\real^d) \to \mathcal P(\real^d)$ is defined by
\[
    \P \nu(v) = \frac{1}{\sqrt{(2\pi)^d \det \Sigma}}\int_{\real^d} \exp \left( - \frac{1}{2} |v - \Psi(u)|_{\Sigma}^2 \right) \nu(u) \, \d u,
\]
and linear operator $\Q\colon \mathcal P(\real^d) \to \mathcal P(\real^d \times \real^{K})$
is determined by
\[
    \Q \ppi (v, y) = \frac{1}{\sqrt{(2\pi)^K \det \Gamma}}
    \exp \left( - \frac{1}{2} \bigl\lvert y - h(v) \bigr\rvert_{\Gamma}^2 \right)
    \ppi(v).
\]
The linear operator~$\P$ determines the evolution of the law of Markov process~$\{v_n\}$,
while the linear operator $\Q$ lifts the law of~$v_n$ to the joint space of state and data $(v_n,y_n).$
Lastly, we denote by~$\B(\placeholder;\yd)\colon \mathcal P(\real^d \times \real^{K}) \to \mathcal P(\real^d)$
the nonlinear operator effecting conditioning of a joint random variable $(v,y)$ on observation $y=\yd$:
\[
  \B(\pi;\yd)(v) = \frac{\pi(v,\yd)}{\int_{\real^d} \pi(u,\yd) \, \d u}.
\]
Having introduced the necessary notation,
we define the probability distribution~$\hmu_{n+1}$ through $v_{n+1}|\Yd_n \sim \hmu_{n+1}$ and observe that
\begin{subequations}
    \label{subequations:true_filter}
    \begin{alignat}{2}
        \hmu_{n+1}&=\P \mmu_n, \qquad && v_{n+1}|\Yd_n \sim \hmu_{n+1}\\
        \pi_{n+1}&=\Q \hmu_{n+1}, \qquad &&(v_{n+1},y_{n+1})|\Yd_n \sim \pi_{n+1}\\
        \mmu_{n+1}&=\B(\pi_{n+1};\yd_{n+1}), \qquad &&{\rm conditioning}\,{\rm on}\, y_{n+1}=\yd_{n+1}.
    \end{alignat}
\end{subequations}
This update formula for $\mu_n$ may be decomposed into a prediction step involving the dynamical model through~$\op P$,
and an analysis step based on Bayes' theorem through the composition $\hmu_{n+1} \mapsto \B(\Q \hmu_{n+1};\yd_{n+1})$.
This iterative interleaving of prediction and analysis steps defines an infinite dimensional dynamical system over the space of probability measures,
numerical approximation of which is a major challenge.

\section{Ensemble Kalman Filter}
\label{section:enkf}

The ensemble Kalman filter was originally introduced by Evensen in \cite{evensen1994sequential}.
We refer to the works~\cite{evensen2009data,evensen2022data} for a detailed review.
In this section,
we first present the basic particle formulation of the algorithm.
Next, we show how this formulation can be derived as an approximation of a mean-field dynamics,
an illuminating perspective fleshed out in~\cite{2022arXiv220911371C}.
Finally, we present stability and error estimates for the mean field filter,
and show how these can be used to study the connection between the finite particle ensemble Kalman filter and the true filter~\cite{carrillo2022ensemble}.

Before presenting the particle formulation of the ensemble Kalman filter, we introduce useful notation:
for a probability measure $\pi \in \mathcal P(\real^d \times \real^{K})$,
we write the mean under~$\pi$ as $\mathcal M(\pi)$ and the covariance under $\pi$ as
\[
    \mathcal C(\pi) =
    \begin{pmatrix}
        \mathcal C^{vv}(\pi) & \mathcal C^{vy}(\pi) \\
        \mathcal C^{vy}(\pi)^\top & \mathcal C^{yy}(\pi)
    \end{pmatrix}.
\]
For $h\colon \real^{d} \to \real^{K}$ we also define ${\mathcal C}^{hh}(\pi)$ as the covariance matrix of the vector $h(v)$,
for $v$ distributed according to the marginal of $\pi$ on the state,
and ${\mathcal C}^{vh}(\pi)$ as the covariance between $v$ and $h(v)$.
An iteration of the ensemble Kalman filter then reads, for $n \in \mathbb{N}$,
\begin{subequations}
    \label{subequations:enkf}
\begin{alignat}{2}
            \widehat v_{n+1}^{(j)} &= \Psi \bigl(v_n^{(j)}\bigr) + \xi_n^{(j)},\\ %\quad v_{n}^{(j)} \sim \muk_{n},\\
            \widehat y_{n+1}^{(j)} &= h(\widehat v_{n+1}^{(j)}) + \eta_{n+1}^{(j)}, \\
            v_{n+1}^{(j)} &= \widehat v_{n+1}^{(j)} + \mathcal C^{vh}\bigl(\pij_{n+1}\bigr) \left( \Gamma + \mathcal C^{hh}\bigl(\pij_{n+1}\bigr) \right)^{-1} \Bigl(\yd_{n+1} - \widehat y_{n+1}^{(j)} \Bigr),\\
\pij_{n+1}&=
\frac{1}{J}\sum_{j=1}^J \delta_{\bigl(\widehat v_{n+1}^{(j)},\widehat y_{n+1}^{(j)}\bigr)}.
\end{alignat}
\end{subequations}
Here $\xi_{n}^{(j)} \sim \normal(0, \Sigma)$ and $\eta_{n}^{(j)} \sim \normal(0, \Gamma)$ are i.i.d.\ with respect to both $n$ and $j$,
and the two sets are independent.
For $n=0$ we choose independent $v_{0}^{(j)} \sim \mmu_0$.
From the particles evolving according to iteration~\eqref{subequations:enkf},
we define the following empirical measure as an approximation of the filtering distribution:
\begin{equation}
    \label{eq:enkf_approximation}
    \muj_{n}= \frac{1}{J}\sum_{j=1}^J \delta_{v_{n}^{(j)}}.
\end{equation}
Since an equal weight is assigned to all the particles,
the ensemble Kalman filter cannot suffer from weight collapse,
which contributes to its popularity among practitioners.
Furthermore, in the setting where the functions~$\Psi,h$ are both affine and $\mu_0$ is Gaussian,
the empirical approximation $\muj_{n}$ converges to $\mmu_n$ at the Monte Carlo rate in the limit as $J \to \infty$,
in an appropriate metric~\cite{le2009large,mandel2011convergence}.
In this setting, the ensemble Kalman filter is an alternative to the classical Kalman filter~\cite{kalman1960new},
the latter being an iterative algorithm for updating sequentially the mean and covariance of the (Gaussian) filtering distribution.

We now present recent research aimed at studying the statistical accuracy of the ensemble Kalman filter beyond the linear Gaussian setting.
To bound the distance between the true filter and its approximation~\eqref{eq:enkf_approximation},
a natural approach is use as a pivot, in a triangle inequality,
the mean field limit of the ensemble Kalman particle system, obtained by letting~$J \to \infty$ in~\eqref{subequations:enkf}.
This limit leads to the following mean field dynamics:
\begin{subequations}
    \label{subequations:enkf_meanfield}
\begin{alignat}{2}
    &\widehat v_{n+1} = \Psi (v_n) + \xi_n,\\
    &\widehat y_{n+1} = h(\widehat v_{n+1}) + \eta_{n+1}, \\
    & v_{n+1} = \widehat v_{n+1} + \mathcal C^{vy}\bigl(\nuk_{n+1}\bigr) \mathcal C^{yy}\bigl(\nuk_{n+1}\bigr)^{-1} \bigl(\yd_{n+1} - \widehat y_{n+1} \bigr),\\
    & (\widehat v_{n+1}, \widehat y_{n+1}) \sim \nuk_{n+1}.
\end{alignat}
\end{subequations}
Once again $\xi_{n} \sim \normal(0, \Sigma)$ and $\eta_{n} \sim \normal(0, \Gamma)$ are i.i.d.\ with respect to $n$,
and the two sets of random variables are independent.
Equations~\eqref{subequations:enkf_meanfield} define a stochastic map from $v_n$ to $v_{n+1}$.
This map is also  nonautonomous, as it depends on the observed data,
and mean field in the sense that $v_{n+1}$ depends not only on~$(\widehat v_{n+1},\widehat y_{n+1})$,
but also on the law~$\nuk_{n+1}$ of this random vector.
In order to write an evolution equation for the law of~$v_{n}$,
which we denote by~$\muk_n$,
it is useful to define $\mapT$ as follows:
\begin{alignat*}{2}
\mapT(\placeholder, \placeholder; \pi,\yd) &\colon
        \real^{d} \times \real^{K} \to \real^d; \\
        (v, y) &\mapsto v + \mathcal C^{vy}(\pi) \, \mathcal C^{yy}(\pi)^{-1} \bigl(\yd - y \bigr).
\end{alignat*}
This map is affine in its first two arguments for any fixed pair $(\pi,\yd)\in \mathcal P(\real^d\times\real^K)\times\real^K$.
We then define the following nonlinear operator on probability measures:
\begin{alignat*}{2}
\T(\pi;\yd) &= \Bigl(\mapT(\placeholder, \placeholder; \pi,\yd)\Bigr)_\sharp {\pi},
\end{alignat*}
where $F_\sharp {\nu}$ denotes the pushforward of a probability measure $\nu$ under a map~$F$.
It is shown in~\cite{2022arXiv220911371C} that the evolution of the probability measure~$\muk_{n}$ is determined by the dynamics~\eqref{eq:kalman_mf} following.
To understand the relationship between $\muk_n$ and $\mmu_n$,
we write the two evolutions in parallel:
\begin{subequations}
\begin{align}
    \label{eq:kalman_mf}
    \muk_{n+1} &= \T(\Q\P\muk_n;\yd_{n+1}),  \qquad \muk_0=\mmu_0 \\
    \label{eq:true_filter}
    \mmu_{n+1} &= \B(\Q\P\mmu_n;\yd_{n+1}).
\end{align}
\end{subequations}
Thus, in order to control the distance between~$\muk_{n+1}$ and~$\mmu_{n+1}$ given a control at iteration~$n$,
it is essential to understand when~$\T \approx \B.$
In fact,
on the set of Gaussian measures $\mathcal G(\R^d \times \R^K)$,
operator~$\T$ performs an exact transport from the joint distribution of state and data to the conditional distribution of state given data~\cite{2022arXiv220911371C};
that is to say~$\T(\pi, y^{\dagger}) = \B(\pi, y^{\dagger})$ for all $\pi \in \mathcal G(\R^d \times \R^K)$ and~$y^{\dagger} \in \real^K$.
This reflects the fact that,
in the linear Gaussian setting,
the ensemble Kalman filter is a particle approximation of the Kalman filter that is asymptotically exact as~$J \to \infty$.

To obtain an error bound for the approximate filtering distribution,
it is useful to define operator~$\G\colon \mathcal{P}_2~\to~\mathcal{G}$,
where $\mathcal P_2$ denotes the space of probability distributions with finite second moment,
mapping an arbitrary probability distribution to the Gaussian distribution with matching first and second moments:
\begin{align*}
    \G\ppi &= \normal\Bigl(\mathcal M(\ppi), \mathcal C(\ppi)\Bigr).
\end{align*}
The probability distribution~$\G \pi$ may also be characterized as the Gaussian distribution~$\pp$ such that Kullback--Leibler divergence~$d_{\rm KL}(\pi\|\pp)$ of~$\pi$ from~$\pp$ is minimized~\cite[Theorem 4.7]{2018arXiv181006191S}.
To describe the main contributions of the works~\cite{carrillo2022ensemble,calvello2024accuracyensemblekalmanfilter},
we also introduce the following weighted variation metric, with weight $g(v) = 1 + |v|^2$:
\[
    \tvdg(\mmu, \mnu) = \sup_{|f| \leq g} \Bigl\lvert \mmu[f] - \mnu[f] \Bigr\rvert,
\]
Here $\mmu[f] = \int f \d \mu$ and similarly for $\mnu[f]$.
For a function $f \colon \real^r \to \real$,
let also $\winf{f}$ denote the following weighted $L^{\infty}$ norm of~$f$,
which is finite if and only if~$f$ grows at most linearly:
\[
    \winf{f} := \esup_{x\in \real^r}  \frac{f(x)}{1 + |x|}.
\]
Finally, for a globally Lipschitz continuous function~$f$,
let $|f|_{C^{0,1}}$ denote its Lipschitz constant.
Armed with this
we proceed to summarize the main results from~\cite{carrillo2022ensemble,calvello2024accuracyensemblekalmanfilter}.
For precise statement and proof of the following theorems,
and more details on the conditions under which they hold,
we refer the reader to these references.
The first result is a stability estimate for the mean field ensemble Kalman filter.
\begin{theorem}
    [Stability: Mean Field Ensemble Kalman Filter]
    \label{theorem:main_theorem}
    Assume the probability measures~$(\muk_n)_{n \in \range{0}{N}}$ and  $(\mu_n)_{n \in \range{0}{N}}$ are obtained from
    the dynamical systems~\eqref{eq:kalman_mf} and~\eqref{eq:true_filter}, respectively.
    Suppose also that functions $\Psi,h$ satisfy $\winf{\Psi}, \winf{h} \leq \kappa$ for some $\kappa > 0$
    and~$|h|_{C^{0,1}} \leq \ell$ for some $\ell > 0$,
    and that the noise covariance matrices~$\Sigma, \Gamma$ are positive definite.
    Then there exists~$C = C(\Sigma, \Gamma, \kappa, \ell, N) > 0$
    such that
    \begin{equation}
        \label{eq:stability}
        \tvdg\bigl(\muk_N, \mu_N\bigr) \leq C  \max_{n \in \range{0}{N-1}} \tvdg( \op Q \op P \mu_n, \op G \op Q \op P \mu_n).
    \end{equation}
\end{theorem}
The right-hand side of~\eqref{eq:stability},
which by this theorem upper bounds the error affecting the mean field ensemble Kalman filter,
is a measure of how close the true filter $\{\mmu_n\}$ is to being Gaussian.
This result was proved in~\cite{carrillo2022ensemble} under uniform boundedness assumptions for~$\Psi,h$,
and then generalized in~\cite{calvello2024accuracyensemblekalmanfilter} to the form presented here to allow linear growth of these functions.
To make this result useful, we must identify concrete settings in which the right-hand side of~\eqref{eq:stability} is indeed small.
In~\cite{carrillo2022ensemble} the setting where~$\Psi, h$ are small perturbation of constant functions is considered.
Here, we present a more general result from the subsequent work~\cite{calvello2024accuracyensemblekalmanfilter},
where small perturbations of affine functions are considered.
With this goal, for functions $f, g$ and $r\geq 0$,
we denote by~$B_{L^{\infty}}\bigl((f,g),r\bigr)$ the~$L^\infty$ ball of radius $r$ centered at $(f, g)$,
on the product space.
\begin{theorem}
    [Error Estimate: Mean Field Ensemble Kalman Filter]
    \label{theorem:error_enkf}
    Assume again that the probability measures~$(\muk_n)_{n \in \range{0}{N}}$ and  $(\mu_n)_{n \in \range{0}{N}}$ are obtained from the dynamical systems~\eqref{eq:kalman_mf} and~\eqref{eq:true_filter}, respectively.
    Let $\Psi_0, h_0$ be affine functions satisfying $\winf{\Psi_0}, \winf{h_0} \leq \kappa$
    and $\lipnorm{h_0} \leq \ell$ for some $\kappa, \ell > 0$,
    and suppose that the noise covariances~$\Sigma, \Gamma$ are positive definite.
    Then there exists a constant $C = C(\Sigma, \Gamma, \kappa, \ell, N) > 0$ such that
    for all~$\varepsilon\in[0,1]$ and $(\Psi,h) \in B_{L^\infty}\bigl((\Psi_0,h_0),\varepsilon \bigr)$
    satisfying $\winf{\Psi}, \winf{h} \leq \kappa$ and $|h|_{C^{0,1}} \leq \ell$,
    it holds that
    \[
        \tvdg(\muk_N, \mu_N) \leq C\varepsilon.
    \]
\end{theorem}

The two preceding theorems concern the mean field limit of the ensemble Kalman filter.
To quantify the error for the finite particle approximation of this limit given in~\eqref{subequations:enkf},
we rely on a triangle inequality,
using a propagation of chaos estimate from~\cite{le2009large} to bound the difference between mean field and finite particle formulations.
This leads to the following bound on the approximation error for the finite particle ensemble Kalman filter.
\begin{theorem}
    [Error Estimate: Finite Particle Ensemble Kalman Filter]
    \label{theorem:finite_enkf_theorem}
    Let $(\mu_n)_{n \in \range{0}{N}}$ denote the true filtering distribution,
    initialized at a Gaussian probability measure $\mu_0\in~\mathcal G(\real^{d})$.
    Denote by $(\muj_n)_{n \in \range{0}{N}}$ its approximation as in~\eqref{eq:enkf_approximation} by the finite particle ensemble Kalman filter~\eqref{subequations:enkf},
    with each particle initialized independently according to~$\mu_0$.
    Then, under the same assumptions as in Theorem~\eqref{theorem:error_enkf},
    with the additional assumption that~$h$ is an affine function,
    there is $C = C(\Sigma, \Gamma, \kappa, \ell, N) > 0$ such that
    for all~$\varepsilon\in[0,1]$ and all~$(\Psi,h) \in B_{L^\infty}\bigl((\Psi_0,h_0),\varepsilon \bigr)$
    satisfying $\winf{\Psi}, \winf{h} \leq \kappa$ and~$|h|_{C^{0,1}} \leq \ell$,
    it holds for any $\varphi$ in an appropriate class of test functions that
    \[
       \left(\mathbb{E}~\Bigl|\muj_N[\varphi] -\mu_N[\varphi]\Bigr|^2\right)^{\frac{1}{2}} \leq C  \left(\frac{1}{\sqrt{J}} + {\varepsilon} \right).
    \]
\end{theorem}

The results presented here are first steps to analyze the accuracy of the ensemble Kalman filter beyond the Gaussian setting.
To conclude, we identify just a few promising research avenues for future work.
We refer to the conclusion sections in~\cite{carrillo2022ensemble,calvello2024accuracyensemblekalmanfilter} for more thorough discussions.
\begin{itemize}
    \item
        First, let us note that in all the theorems presented here,
        the constant prefactor~$C(\Sigma, \Gamma, \kappa, \ell, N)$ on the right-hand side grows exponentially with the number of steps~$N$.
        This is partly because the proof of~Theorem~\ref{theorem:main_theorem} follows a classical ``consistency plus stability implies convergence" approach,
        which ultimately relies on a Gr\"onwall inequality.
        It would be of great interest to improve the results so that they remain useful in the large~$N$ setting,
        which will require to understand long-time stability properties of the true filter.

    \item
        The ensemble Kalman filter may be used, in a modified form,
        as a derivative-free method to solve inverse problems.
        This application of the filter has received sustained attention from the mathematical community in recent years,
        see for example~\cite{MR3041539,MR4059375,MR4482475,calvello2024accuracyensemblekalmanfilter}.
        It would be interesting to investigate the extent to which our results can be generalized to this setting.

    \item
        We focused in the highlighted works exclusively on the discrete-time setting,
        but a continuous-time analysis would also be of interest.
        Propagation of chaos results for the Kalman--Bucy filter, such as those from paper~\cite{MR3784489},
        may prove important in this context.
\end{itemize}

{\small
\bibliographystyle{abbrv}
\bibliography{template}
}

\end{document}